\documentclass[]{interact}

\usepackage{epstopdf}% To incorporate .eps illustrations using PDFLaTeX, etc.
\usepackage[caption=false]{subfig}% Support for small, `sub' figures and tables

\usepackage[numbers,sort&compress]{natbib}% Citation support using natbib.sty
\bibpunct[, ]{[}{]}{,}{n}{,}{,}% Citation support using natbib.sty
\renewcommand\bibfont{\fontsize{10}{12}\selectfont}% Bibliography support using natbib.sty

\theoremstyle{plain}% Theorem-like structures provided by amsthm.sty

\theoremstyle{definition}

\theoremstyle{remark}

\begin{document}

\articletype{RESEARCH PAPER}% Specify the article type or omit as appropriate

\title{Multi-Gaussian random variables}

\author{
\name{O.~Korotkova\thanks{CONTACT O.~Korotkova Email: korotkova@physics.miami.edu}} 
\affil{Department of Physics, University of Miami, 1320 Campo Sano Dr., Coral Gables, FL 33146}}

\maketitle

\begin{abstract}
A generalization of the classic Gaussian random variable to the family of Multi-Gaussian (MG) random variables characterized by shape parameter $M>0$, in addition to the mean and the standard deviation, is introduced. The probability density function of the MG family members is the alternating series of the Gaussian functions with the suitably chosen heights and widths. In particular, for the integer values of $M$ the series has finite number of terms and leads to flattened profiles, while reducing to classic Gaussian density for $M=1$. For non-integer, positive values of $M$ a convergent infinite series of Gaussian functions is obtained that can be truncated in practical problems. While for all $M>1$ the MG PDF has flattened profiles, for $0<M<1$ it leads to cusped profiles. Moreover, the multivariate extension of the MG random variable is obtained and the Log-Multi-Gaussian (LMG) random variable is introduced. 
\end{abstract}


\section{Introduction}

The most famous Probability Density Function (PDF) of a continuous random variable - Gaussian - stemming from the early works of de Moivre \cite{Moivre} and Gauss \cite{Gauss} can be generalized in a number of ways for inclusion of desired shape details such as flattening, skewing, splitting, etc. One well-known generalization was proposed by Subbotin \cite{Subbotin} (see also \cite{Levy}) who extended the PDF curve to flatter or sharper versions by varying the power law of the exponential function, hence the family is sometimes termed exponential power or super-Gaussian. The Subbotin's PDF was later rescaled by Lunetta \cite{Lunetta} and the resulting family has become well explored (c.f. \cite{Nadarajah}). While originally they were used for analysis of the astrophysical data, the super-Gaussian random variables are currently employed for characterization of a wide range of statistical phenomena: in big data analysis in general \cite{Lah} and, in particular, in finance \cite{Mantegna}, \cite{McCauley}, genetics \cite{Liang} and scientific impact assessment \cite{Redner}, to name a few. However, this seemingly transparent generalization often rely on the use of special functions, such as a hypergeometric function, as is the case for evaluation of its characteristic function \cite{CFsg} (see also \cite{Maturi}). 

In this paper we introduce a novel family of continuous random variables that serves a similar purpose as the super-Gaussian family, i.e., it reshapes the Gaussian distribution to flat-top or cusp-top versions, depending on the value of the shape parameter. However, the main advantage of our family over the super-Gaussians stems from the fact that the statistical properties of its members can be expressed as the series of those for the Gaussian random variables, with very simple expressions defining their heights  and widths. Moreover, in case when the shape parameter is an integer the flat-top distributions can be formed by a finite number of terms in the series. As we show, the ability to represent a PDF of the new random variable as a linear combination of Gaussian contributions leads to unprecedented tractability in derivation of a number of its characteristics. For any value of the shape parameter we will term our new family of random variables \textit{Multi-Gaussian} (MG), not to be confused with the well-known multivariate Gaussian. The multi-Gaussian functions of various dimensions have been previously used in optics for modeling of beam intensity profiles \cite{MG0}, aperure shapes \cite{MG1}, scattering potentials \cite{MG2},  \cite{MG3} and various correlation and coherence functions \cite{MG4}-\cite{MG7}.  

Starting from the finite series case (integer values of shape parameter) we first build in some intuition for the new random variable. In particular, after introducing the PDF, we provide calculations of its Cumulative Distribution Function (CDF), Moment Generating Function (MGF), Characteristic Function (CF) and the Cumulant Generating Function (CGF). Then we derive the general expressions for the moments and the cumulants of any order and find explicit expressions for the first four members of each sequence. 

Next, we introduce the Log-Multi-Gaussian (LMG) random variable on assuming that its logarithm is MG-distributed. Such derived distribution can be regarded as an extension of the classic Log-Normal distribution \cite{McAlister} to flat-topped profiles. The LMG distribution has an analog derived from the super-Gaussian family \cite{Vianelli1982}, \cite{Vianelli1983}, however, as we show, the calculations of major statistical characteristics of the LMG random variables can be obtained almost effortlessly, as the linear combinations of the well-known results for the Log-Normal variables. 

We also discuss the natural extension of the univariate MG random variable to the multivariate domain and, in particular, the bivariate domain. Such variables redice to classic multivariate/bivariate Gaussians if the shape parameter $M=1$. 

Finally, we show how to generalize all the aforementioned results to the MG distributions with shape parameter $M$ taking on any positive values, not necessarily integers, while discussing numerical examples for the special case when $M$ is a reciprocal of an integer, leading to formation of various cusped distributions. The features of the LMG and bivariate MG random variables with any positive $M$ are also briefly outlined and the corresponding numerical results are provided.

\section{Multi-Gaussian distribution with integer shape index $M$}

\subsection{Probability Density Function and Cumulative Distribution Function}

Let us begin by recalling that the Gaussian PDF of a continuous real random variable $X$ with mean $\mu$ (location parameter) and standard deviation $\sigma$  (scale parameter), i.e. $X \sim \mathcal{N}(\mu,\sigma)$, has form
\begin{equation}\label{G}
p_X^{(G)}(x)=\frac{1}{\sqrt{2\pi}\sigma}\exp\left[ - \frac{(x-\mu)^2}{2\sigma^2}\right].
\end{equation}
Consider now a function 
\begin{equation}\label{fx}
f(x)=\frac{1}{\sqrt{2\pi}\sigma}\left[1-\left( 1-\exp\left[ - \frac{(x-\mu)^2}{2\sigma^2}\right] \right)^M \right],
\end{equation} 
where $M$ is a positive real number. For $M=1$ $f(x)$ reduces to the Gaussian function in Eq. \eqref{G}, while for $M>1$ and $0<M<1$ it describes flat-top and cusped distributions, respectively. Hence, we may refer to $M$ as a shape parameter.

Let us first discuss the case when $M$ is an integer, $M=1,2,3,...$ On using the binomial theorem 
\begin{equation}
(u+v)^M=\sum\limits_{m=0}^{M} \binom{M}{m}u^{M-m}v^m, \quad \binom{M}{m}=\frac{M!}{m!(M-m)!}, 
\end{equation}  
with  $u=1$, $v=-\exp\left[ - \frac{(x-\mu)^2}{2\sigma^2}\right]$ we arrive at the finite series
\begin{equation}
f(x)=\frac{1}{\sqrt{2\pi}\sigma} \sum\limits_{m=1}^{M} \binom{M}{m} (-1)^{m-1}\exp\left[ - \frac{(x-\mu)^2}{2\sigma_m^2}\right],
\end{equation}
with the standard deviation of the $m$-th term in the series of the form 
\begin{equation}
\sigma_m=\sigma/ \sqrt{m},
\end{equation}
and the same mean $\mu$ for all terms. Let us now introduce a PDF as 
\begin{equation}\label{pxnorm}
p_X^{(MG)}(x)=f(x) / \int\limits_{-\infty}^{\infty}f(x)dx, 
\end{equation}
ensuring that $\int\limits_{-\infty}^{\infty}p_X^{(MG)}(x)dx=1$. On finding that 
\begin{equation}
C_0(M)=\int\limits_{-\infty}^{\infty}f(x)dx=\sum\limits_{m=1}^{M} \binom{M}{m} \frac{(-1)^{m-1}}{\sqrt{m}},
\end{equation}
we finally obtain the \textit{Multi-Gaussian} PDF:
\begin{equation}\label{MGpdf}
p_X^{(MG)}(x)=\frac{1}{C_0(M)\sqrt{2\pi}\sigma} \sum\limits_{m=1}^{M} \binom{M}{m} (-1)^{m-1}\exp\left[ - \frac{(x-\mu)^2}{2\sigma_m^2}\right].
\end{equation}
We may also say that $X \sim \mathcal{N}(\mu,\sigma,M)$ as a generalization of normal variable to any values of shape parameter $M$.

\begin{figure}[h!]
\centering\includegraphics[width=14cm]{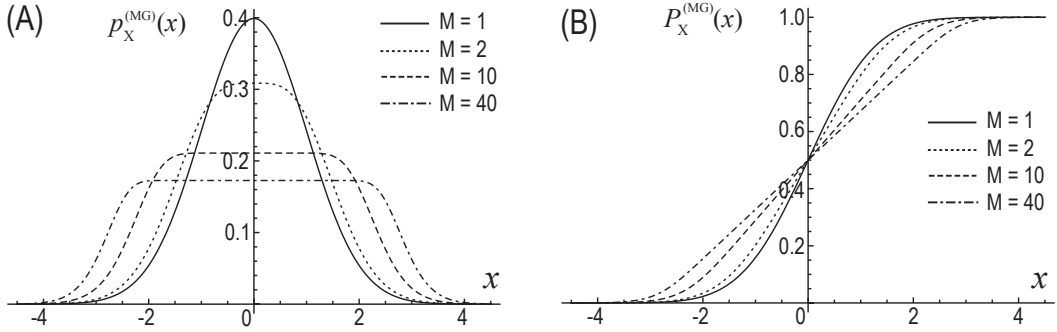}
\caption{Multi-Gaussian PDF (A) and CDF (B) with $\mu=0$, $\sigma=1$. \label{figMG}}
\end{figure}

\begin{figure}[h!]
\centering\includegraphics[width=14cm]{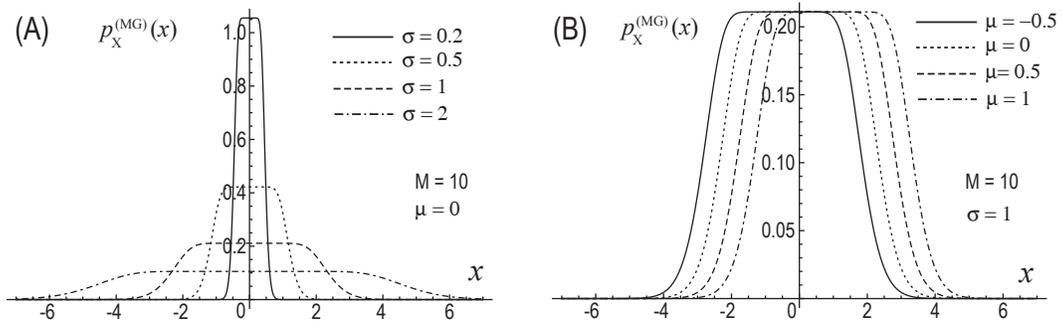}
\caption{Multi-Gaussian PDF with $M=10$ and various values of $\mu$ and $\sigma$.\label{figMGsigmamu}}
\end{figure}

The Cumulative Distribution Function (CDF) of the Multi-Gaussian random variable in Eq. \eqref{MGpdf} can then be readily calculated from its definition and by changing the order of summation and integration: 
\begin{equation}
\begin{split}
P_X^{(MG)}(x)&=\int\limits_{-\infty}^{x} p_X^{(MG)}(s)ds\\&
= \frac{1}{2C_0(M)} \sum\limits_{m=1}^{M} \binom{M}{m}\frac{ (-1)^{m-1}}{\sqrt{m}}\left[1+Erf\left(  \frac{x-\mu}{\sqrt{2}\sigma_m}\right) \right],
\end{split}
\end{equation}
where $Erf$ stands for the error function. Figures \ref{figMG} (A) and (B) show the PDF and the CDF of a MG random variable for $M=1,2,10$ and $40$, calculated from Eqs. (8) and (9), respectively. Figure \ref{figMGsigmamu} shows the PDF for $M=10$ and different values of $\mu$ and $\sigma$.

\subsection{Moment Generating Function, Characteristic Function and Moments}

The Moment Generating Function (MGF) of the MG random variable in Eq. \eqref{MGpdf} can be also directly evaluated from its definition: 
\begin{equation}\label{MGF}
\begin{split}
M_X(t)^{(MG)}&=\int\limits_{-\infty}^{\infty} \exp[x t]p_X(x)dx \\&
= \frac{1}{C_0(M)}\exp[\mu t] \sum\limits_{m=1}^{M} \binom{M}{m}\frac{ (-1)^{m-1}}{\sqrt{m}} \exp\left[\frac{\sigma_m^2 t^2}{2} \right].
\end{split}
\end{equation} 

The Characteristic Function (CF) of the MG distribution can be readily obtained either from definition or  its relation with the MGF:
\begin{equation}
\begin{split}
\phi_X^{(MG)}(\omega)&=\int\limits_{-\infty}^{\infty} \exp[i\omega x]p_X^{(MG)}(x)dx =M_X^{(MG)}(i\omega) \\&
= \frac{1}{C_0(M)}\exp[i\omega \mu] \sum\limits_{m=1}^{M} \binom{M}{m}\frac{ (-1)^{m-1}}{m^n\sqrt{m}} \exp\left[\frac{-\sigma_m^2 \omega^2}{2}\right].
\end{split}
\end{equation}

The statistical moment of order $k$ can be either evaluated from the MGF in Eq. \eqref{MGF} via expression 
\begin{equation}
\mu_k^{(MG)}=\frac{d^k}{dt^k}[M_X(t)]_{t=0},
\end{equation}
or directly via the PDF function in Eq. \eqref{MGpdf}. Following the latter path we find that
\begin{equation}\label{moments}
\begin{split}
\mu_k^{(MG)}&=\int\limits_{-\infty}^{\infty} x^k p_X^{(MG)}(x)dx \\&
=\frac{1}{C_0(M)} \sum\limits_{m=1}^{M} \binom{M}{m}\frac{ (-1)^{m-1}}{\sqrt{m}} \int\limits_{-\infty}^{\infty} \frac{x^k}{\sqrt{2\pi}\sigma_m}\exp\left[-\frac{(x-\mu)^2}{2\sigma_m^2} \right]dx,
\end{split}
\end{equation}
where we recognize that each integral term is the $k$-th moment of the Gaussian distribution [see Eq. (1)] with mean $\mu$ and standard deviation $\sigma_m$:
\begin{equation}
\mu_k^{(G)}=\sigma^2(-i\sqrt{2})^kU\left( \frac{k}{2},\frac{1}{2},-\frac{\mu^2}{2\sigma_m^2}\right),
\end{equation}
where $U(a,b,z)$ is the confluent Hypergeometric function. In particular, on substituting the first four moments of the Gaussian distribution from Eq. (14) into Eq. \eqref{moments} we find at once that
\begin{equation}
\begin{split}
&\mu_1^{(MG)}=\mu, \\&
\mu_2^{(MG)}=\mu^2+\sigma^2 \xi_1(M), \\&
\mu_3^{(MG)}=\mu^3+3\mu \sigma^2 \xi_1(M), \\&
\mu_4^{(MG)}=\mu^4+6\mu^2\xi_1(M) +3\sigma^4\xi_2(M).
\end{split}
\end{equation}
Here parameters $\xi_n(M)$ are defined via ratios
\begin{equation}
\xi_n(M)=\frac{C_n(M)}{C_0(M)}, \quad n=0,1,2,...
\end{equation}
where
\begin{equation}\label{Cn}
C_n(M)=\sum\limits_{m=1}^{M} \binom{M}{m}\frac{ (-1)^{m-1}}{m^n\sqrt{m}},
\end{equation}
and, in particular, in agreement with $C_0(M)$ defined in Eq. (7).  Also, $\xi_n(1)=1$ for any $n=1,2,3,...$, hence sequence (15) reduces to that in Eq. (14).

\subsection{Cumulant Generating Function and Cumulants}
 
The Cumulant Generating Function (CGF) of the MG random variable has form 
\begin{equation}
\begin{split}
K_X^{(MG)}(h)&=\ln \left[M_X^{(MG)}(h) \right] \\&
=h\mu+\ln\left[ \frac{1}{C_0(M)}\sum\limits_{m=1}^{M} \binom{M}{m}\frac{ (-1)^{m-1}}{\sqrt{m}} \exp\left(\frac{h^2 \sigma_m^2}{2}\right)\right],
\end{split}
\end{equation}
as implied by Eq. (10). On expanding the exponential function in Eq. (18) in Taylor series, interchanging the order of two summations and recognizing coefficients $\xi_n(M)$ in the sum inside we get
\begin{equation}
K_X^{(MG)}(h)=h\mu+\ln\left[ \sum\limits_{n=0}^{\infty} \frac{h^{2n} \sigma^{2n}}{2^n n!} \xi_n(M)\right].
\end{equation}
Further, expanding the logariphmic function in Eq. (19) in the Taylor series we arrive at the double power series
\begin{equation}\label{CGF}
K_X^{(MG)}(h)=\mu h+\sum\limits_{p=1}^{\infty}\frac{(-1)^{p+1}}{p}\left[ \sum\limits_{n=0}^{\infty}  \frac{h^{2n}\sigma^{2n}}{2^n n!}  \xi_n(M) \right]^p.
\end{equation} 
The cumulants can be found as coefficients $\kappa_n^{(MG)}$ of powers of $h$ in expansion (\ref{CGF}):
\begin{equation}
K_X^{(MG)}(h)=\sum\limits_{p=1}^{\infty} k_p^{(MG)} \frac{h^p}{p!}.
\end{equation}
In particular, the first four cumulants are:
\begin{equation}
\begin{split}
&\kappa_1^{(MG)}=\mu, \\&
\kappa_2^{(MG)}=\sigma^2\xi_1(M), \\&
\kappa_3^{(MG)}=0, \\&
\kappa_4^{(MG)}=3\sigma^4[\xi_2(M)-\xi_1^2(M)].
\end{split}
\end{equation}
and, all cumulants of odd orders higher than three are also trivial. 

\section{Log-Multi-Gaussian (LMG) distribution}

\begin{figure}[h!]
\centering\includegraphics[width=14cm]{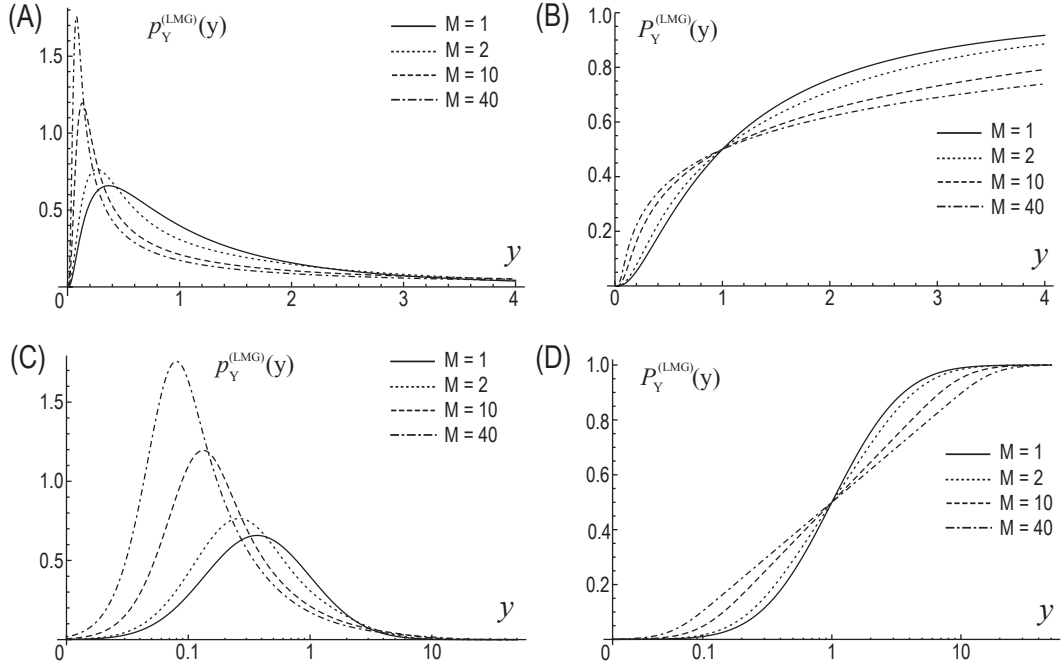}
\caption{Log-Multi-Gaussian PDF (A), (C) and CDF (B), (D) with $\mu=0$, $\sigma=1$. (C) and (D) are the same as (A) and (B) but on log-linear scale. \label{figLMGM}}
\end{figure}

\begin{figure}[h!]
\centering\includegraphics[width=14cm]{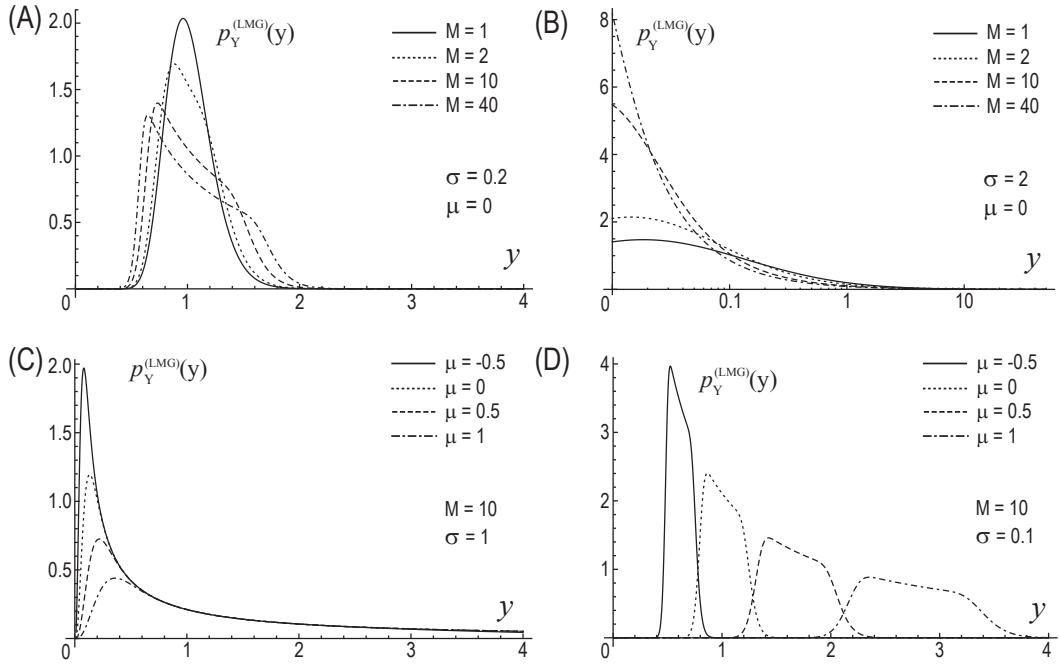}
\caption{Log-Multi-Gaussian PDF with different values of $\sigma$ and $\mu$. \label{figLMGsigmamu}}
\end{figure}

Let random variable $X$  be the MG-distributed and $Y$ is such that $\ln Y=X$. Then, if we let $Y=\exp[X]=g(X)$ and $X=\ln[Y]=g^{-1}(Y)$, then the PDF for variable Y takes form
\begin{equation}
\begin{split}\label{logMGpdf}
p_Y(y)&=p_X[g^{-1}(y)] \left| \frac{d g^{-1}(y)}{dy}\right| \\&
=p_X[\ln y]\frac{1}{|y|}, \quad y>0.
\end{split}
\end{equation}
Substitution of the MG PDF from Eq. \eqref{MGpdf} into Eq. \eqref{logMGpdf} leads to PDF
\begin{equation}
p^{(LMG)}_Y(y)=\frac{1}{C_0(M)\sigma \sqrt{2\pi}|y|}\sum\limits_{m=1}^M \binom{M}{m} (-1)^{m-1}
\exp\left[-\frac{(\ln y-\mu)^2}{2\sigma_m^2} \right],
\end{equation}
which we will term \textit{Log-Multi-Gaussian (LMG)}. We may say that $Y \sim \mathcal{LN}(\mu, \sigma, M)$. For $M=1$ it reduces to the classic  Log-Normal distribution \cite{McAlister}. The CDF of the LMG distribution can be readily found from relation
\begin{equation}
\begin{split}
P_Y^{(LMG)}(y)&=P_X^{(MG)}(\ln y)\\&
=\frac{1}{2C_0(M)}\sum\limits_{m=1}^M \binom{M}{m}\frac{(-1)^{m-1}}{\sqrt{m}}\left[1+Erf\left[ \sqrt{\frac{m}{2}}\left(\frac{\ln y -\mu}{\sigma}\right)\right] \right], \quad y>0.
\end{split}
\end{equation}
The MGF of the LMG distribution diverges but the moments can be found from definition:
\begin{equation}
\begin{split}
\mu_k^{(LMG)}&=\int \limits_{-\infty}^{\infty}y^k p_Y^{(LMG)}(y)dy \\&
=\frac{1}{C_0(M)} \sum\limits_{m=1}^M \binom{M}{m} \frac{(-1)^{m-1}}{\sqrt{m}} \left\{ \frac{1}{\sqrt{2\pi} \sigma_m} \int \limits_{-\infty}^{\infty}y^{k-1} \exp\left[-\frac{(\ln y-\mu)^2}{2\sigma_m^2} \right] dy \right\}.
\end{split}
\end{equation}
The expression in the curly bracket is the $k$-th moment of the LG distribution with mean $\mu$ and variance $\sigma_m$:
\begin{equation}
\begin{split}
\mu_k^{(LN)}&=\frac{1}{\sqrt{2\pi} \sigma_m} \int \limits_{-\infty}^{\infty}y^{k-1} \exp\left[-\frac{(\ln y-\mu)^2}{2\sigma_m^2} \right] dy \\&
=\exp\left[ \frac{k(2\mu+k\sigma_m^2)}{2}\right] 
\end{split}
\end{equation}
Thus,
\begin{equation}
\mu_k^{(LMG)}=\frac{1}{C_0(M)} \sum\limits_{m=1}^M \binom{M}{m} \frac{(-1)^{m-1}}{\sqrt{m}} \exp\left[ \frac{k(2\mu+k\sigma_m^2)}{2}\right]. 
\end{equation}
In particular, the first four moments are:
\begin{equation}
\begin{split}
&\mu_1^{(LMG)}=\frac{\exp[\mu]}{C_0(M)}\sum\limits_{m=1}^M \binom{M}{m} \frac{(-1)^{m-1}}{\sqrt{m}}\exp\left[\frac{\sigma_m^2}{2}\right],  \\&
\mu_2^{(LMG)}=\frac{\exp[2\mu]}{C_0(M)}\sum\limits_{m=1}^M \binom{M}{m} \frac{(-1)^{m-1}}{\sqrt{m}}\exp\left[2\sigma_m^2\right], \\&
\mu_3^{(LMG)}=\frac{\exp[3\mu]}{C_0(M)}\sum\limits_{m=1}^M \binom{M}{m} \frac{(-1)^{m-1}}{\sqrt{m}}\exp\left[\frac{9\sigma_m^2}{2}\right],  \\&
\mu_4^{(LMG)}=\frac{\exp[4\mu]}{C_0(M)}\sum\limits_{m=1}^M \binom{M}{m} \frac{(-1)^{m-1}}{\sqrt{m}}\exp\left[8\sigma_m^2\right].
\end{split}
\end{equation}
Figure \ref{figLMGM} shows the PDF and the CDF of a LMG random variable with $\mu=0$ and $\sigma=1$,  for several values of index  $M$. As $M$ increases the PDF profiles become sharper with maxima occuring at smaller values of $y$. Figure \ref{figLMGsigmamu} illustrates the PDF of a LMG random variable with $M=10$ but different values of $\sigma$ and $\mu$. 

\section{Multivariate Multi-Gaussian distribution}

Let $\textbf{X}=(X_1, X_2, ..., X_N)^T \in \mathbb{R}$ be an $N$-dimensional vector of real random variables, $T$ standing for transpose. Also let the $N$-dimensional vector of their mean values be $\pmb{\mu}=(\mu_1,\mu_2, ... \mu_N)^T$ and their $N\times N$ covariance matrix be $\pmb{\Sigma}_m=(Cov[X_i, X_j], 1\leq i, j\leq N)/m$. Then the PDF of the multivariate MG random vector can be defined by straightforward generalization of the multivariate Gaussian PDF as 
\begin{equation}
P_{\pmb{X}}^{(MG)}(x_1,x_2,...,x_N)=\frac{\sum\limits_{m=1}^M \binom{M}{m} (-1)^{m-1}\exp\left[-\frac{1}{2}(\pmb{x}-\pmb{\mu})^T \pmb{\Sigma}^{-1}_m(\pmb{x}-\pmb{\mu}) \right]}{C_0(M)\sqrt{(2\pi)^n}\sqrt{ det[\pmb{\Sigma}_m]}},
\end{equation}
where $det$ stands for determinant of a matrix and power $-1$ denotes matrix inverse. 

\begin{figure}[h!]
\centering\includegraphics[width=12cm]{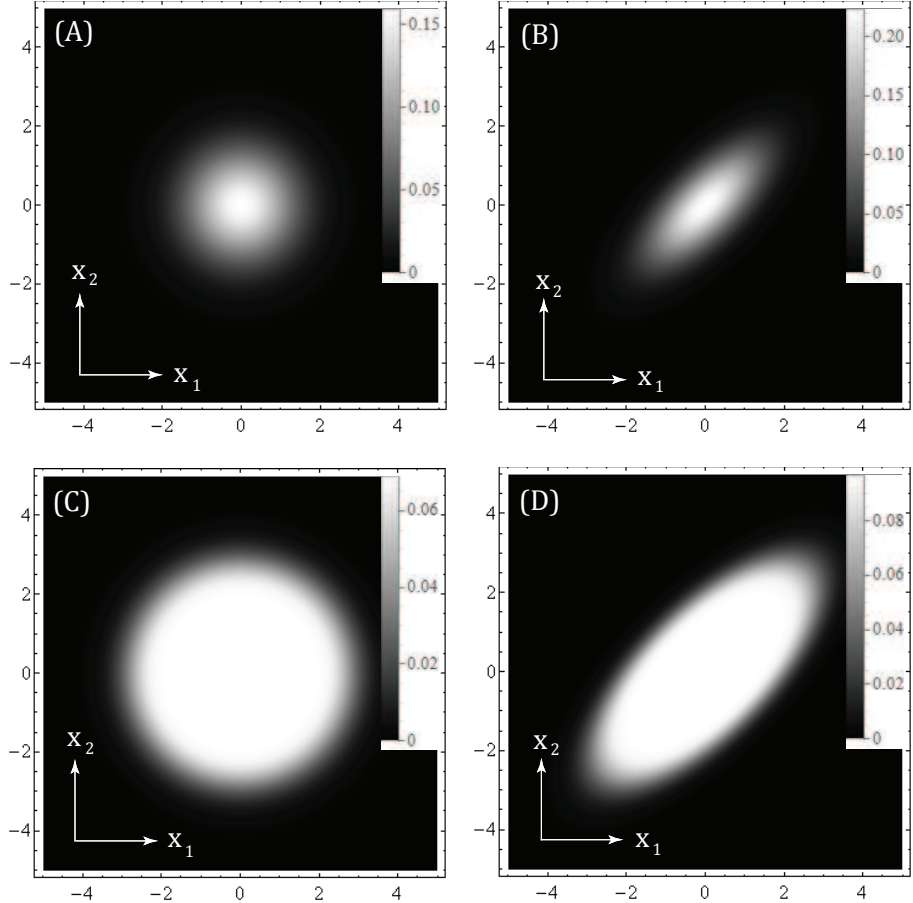}
\caption{Bivariate Multi-Gaussian random variable with $\mu_1=\mu_2=0$, $\sigma_1=\sigma_2=1$. (A) $M=1$, $\rho=0$; (B)  $M=1$, $\rho=0.7$;  (C) $M=40$, $\rho=0$; (D)  $M=40$, $\rho=0.7$. \label{MV}}
\end{figure}

In particular, in the bivariate MG case, $N=2$, one has
\begin{equation}
P_{X_1X_2}^{(MG)}(x_1,x_2)=\frac{\sum\limits_{m=1}^M \binom{M}{m} (-1)^{m-1}\exp\left[-\frac{z m}{2(1-\rho^2)} \right]}{C_0(M) 2\pi \sigma_1\sigma_2 \sqrt{1-\rho^2}},
\end{equation}
with
\begin{equation}
z=\frac{(x_1-\mu_1)^2}{\sigma_{1}^2}-\frac{2\rho(x_1-\mu_1)(x_2-\mu_2)}{\sigma_{1}\sigma_{2}}+\frac{(x_2-\mu_2)^2}{\sigma_{2}^2},
\end{equation}
where we have used 
\begin{equation}
\pmb{\mu}=\begin{pmatrix}
    \mu_{1}  \\
    \mu_{2} 
  \end{pmatrix}, \quad 
\pmb{\Sigma}_m=\frac{1}{m}\begin{pmatrix}
    \sigma_{1}^2 & \rho\sigma_{1}\sigma_{2}  \\
    \rho\sigma_{1}\sigma_{2} & \sigma_{2}^2
  \end{pmatrix}.
\end{equation}
Figure \ref{MV} presents the bivariate MG PDF for several values of $\rho$ and $M$.

\section{Generalization to any $M>0$}

Let us now return to Eq. (2) assuming that $M$ is any positive number, and not necessarily an integer. Then fractional binomial theorem can be used instead of the usual one, hence Eq. (3) becomes:
\begin{equation}\label{uvgen}
(u+v)^M=\sum\limits_{m=0}^{M} \binom{M}{m}u^{M-m}v^m, \quad \binom{M}{m}=\frac{(M)_m}{m!},
\end{equation}  
where 
\begin{equation}\label{poch}
(M)_m=M(M-1)...(M-m+1)
\end{equation}
is the Pochhammer symbol. Let us then set $u=1$, $v=-\exp\left[ - \frac{(x-\mu)^2}{2\sigma^2}\right]$ and express function $f(x)$ in Eq. \eqref{fx} via the generalized binomial series:
\begin{equation}\label{fxgen}
f(x)=\sum\limits_{m=1}^{\infty} \frac{(M)_{m}}{m!}(-1)^{m-1}\exp\left[-\frac{(x-\mu)^2}{2\sigma_m^2} \right].
\end{equation} 
The MG PDF can be obtained as the normalized version of $f(x)$ in Eq. \eqref{fxgen} [see Eq. \eqref{pxnorm}]:
\begin{equation}
p_X^{(MG)}(x)=\frac{1}{C_0(M)\sqrt{2 \pi}\sigma}\sum\limits_{m=1}^{\infty} \frac{(M)_{m}}{m!}(-1)^{m-1}\exp\left[-\frac{(x-\mu)^2}{2\sigma_m^2} \right].
\end{equation}
with normalization
\begin{equation}
C_0(M)=\sum\limits_{m=1}^{\infty} \frac{(M)_{m}}{m!\sqrt{m}}(-1)^{m-1}.
\end{equation}
The infinite series (37) converges as $m \rightarrow \infty$ and hence in practical applications can be suitably truncated. All the caclulations relating to the CDF, MGF, CF, KGF, moments, cumulants, etc. can be carried out in the same manner as done for the integer values of $M$ but with binomial coefficients being of form as in Eqs. \eqref{uvgen} and \eqref{poch}. For instance, the CDF takes form
\begin{equation}
P_X^{(MG)}(x)=\frac{1}{2C_0(M)}\sum\limits_{m=1}^{\infty} \frac{(M)_{m}}{m!\sqrt{m}}(-1)^{m-1} \left[ 1+Erf\left( \frac{x-\mu}{\sqrt{2}\sigma_m}\right)\right].
\end{equation}
\begin{figure}[h!]
\centering\includegraphics[width=14cm]{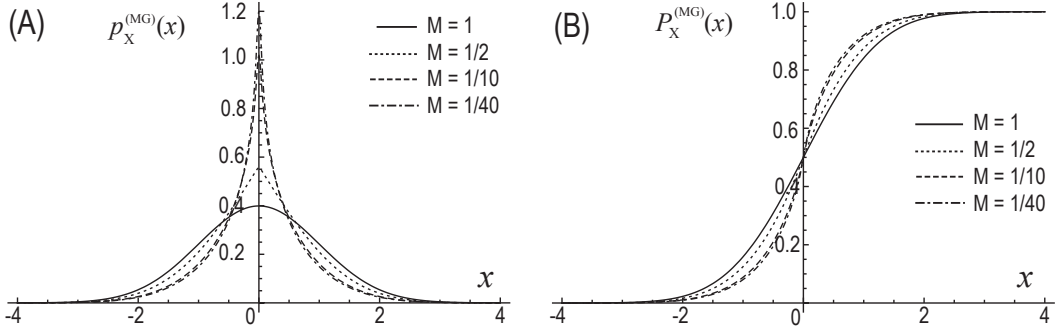}\label{anyMpdfcdf}
\caption{Multi-Gaussian PDF for fractional $M$, with truncation index of $2000$. (A) and CDF (B) with $\mu=0$, $\sigma=1$.}
\end{figure} 

Figure 6 presents the PDF given by Eq. (37) and the CDF given by Eq. (39) of the Multi-Gaussian random variable for several rational values of $M$. 

Likewise, the coefficients $\xi_n(M)$ appearing in calculations of the moments and the cumulants must use more general expression for $C_n(M)$: 
\begin{equation}
C_n(M)=\sum\limits_{m=1}^{\infty} \frac{(M)_m}{m!}\frac{ (-1)^{m-1}}{m^n\sqrt{m}},
\end{equation}
instead of that given in Eq. \eqref{Cn}. 

In particular, the expression for the LMG PDF now takes form 
\begin{equation}
p^{(LMG)}_Y(y)=\frac{1}{C_0(M)\sigma \sqrt{2\pi}|y|}\sum\limits_{m=1}^{\infty} \frac{(M)_{m}}{m!}
\exp\left[-\frac{(\ln y-\mu)^2}{2\sigma_m^2} \right],
\end{equation}
where $C_0(M)$ is given in Eq. (38). Also, the bivariate MG PDF becomes:
\begin{equation}
P_{X_1X_2}^{(MG)}(x_1,x_2)=\frac{1}{C_0(M) 2\pi \sigma_1\sigma_2 \sqrt{1-\rho^2}} \sum\limits_{m=1}^{\infty}  \frac{(M)_{m}}{m!} (-1)^{m-1}\exp\left[-\frac{z m}{2(1-\rho^2)} \right],
\end{equation}
where $C_0(M)$ is given in Eq. (38) but with $z$ being defined as before in Eq. (32). 

Figure 7 present the examples of LMG PDF from Eq. (41) and Fig. 8 shows typical bivariate MG distribution from Eq. (42), both for rational values of index $M$.

\section{Summary}

We have introduced a new family of continuous, real random variables whose PDF represent flattened and cusped deviations from a Gaussian random variable depending on a single shape parameter taking positive values. Gaussian random variables are a particular case of our family when shape parameter takes on value one. While general analytic form of the new family can itself be used it appears possible to express it as a linear combination of Gaussian functions with alternating signs, the same mean and the monotonically decreasing standard deviations. Hence the suggested name of the random variable: Multi-Gaussian. It was shown that for the integer values of the shape parameter the series of Gaussian functions has the finite number of terms while for non-integer positive values it becomes infinite but remains convergent. Due to this feature the calculations of the statistical properties of the Multi-Gaussian random variables essentially reduce to algebraic operations over the well-known properties for a Gaussian random variable. 

We have also introduced the Log-Multi-Gaussian random variable and carried out the extension from univariate Multi-Gaussian random variable to its multivariate counterpart. Other straightforward extensions can be readily made, for instance, relating to complex Multi-Gaussian random variables. The Multi-Gaussian family of random variables is invisioned to be useful in the same applications as the Subbotin's family (exponential-power or super-Gaussian family) but in situations where simple, closed-form analytical results are of importance.

\begin{figure}[h!]
\centering\includegraphics[width=12cm]{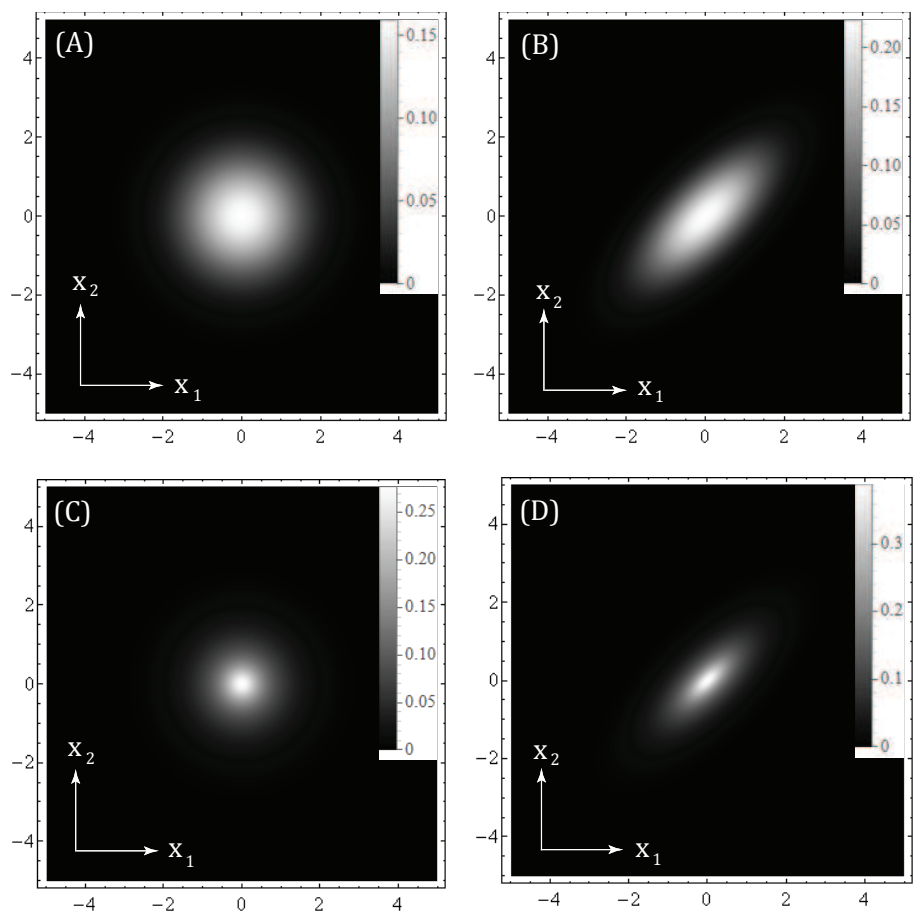}
\caption{Bivariate LMG random variable with $\mu_1=\mu_2=0$, $\sigma_1=\sigma_2=1$. (A) $M=1$, $\rho=0$; (B)  $M=1$, $\rho=0.7$;  (C) $M=1/40$, $\rho=0$; (D)  $M=1/40$, $\rho=0.7$. \label{MVanyM}}
\end{figure}
\begin{figure}[h!]
\centering\includegraphics[width=14cm]{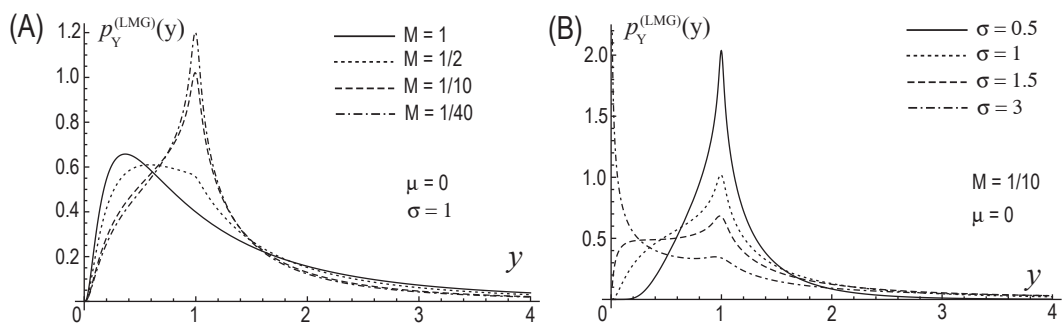}\label{LMGsigmamuAnyM}
\caption{LMG PDF for rational $M$, with truncation index of $2000$, for various values of $M$ $\mu$, and $\sigma$.}
\end{figure}

\end{document}